\documentclass[12pt]{amsart}

\usepackage{amsmath,amsfonts,amssymb}

\def\wtil{\widetilde}
\def\es{\emptyset}
\def\diam{{\rm diam}}

\def\Lk{{\mathcal L}}

\def\lam{\lambda}

\newcommand{\La}{\Lambda}
\newcommand{\be}{\beta}

\newcommand{\e}{\varepsilon}

\newcommand{\BR}{\mathbb{R}}
\newcommand{\BC}{\mathbb{C}}

\newcommand{\BZ}{\mathbb{Z}}

\renewcommand{\limsup}{\varlimsup}
\renewcommand{\liminf}{\varliminf}

\renewcommand{\Re}{{\mathrm{Re}}\,}

\newtheorem{lemma}{Lemma}[section]

\newtheorem{thm}[lemma]{Theorem}
\newtheorem{cor}[lemma]{Corollary}
\theoremstyle{definition}
\newtheorem{Def}[lemma]{Definition}
\newtheorem{example}[lemma]{Example}

\theoremstyle{remark}
\newtheorem{rmk}[lemma]{Remark}
\newtheorem{rmks}[lemma]{Remarks}

\numberwithin{equation}{section} \numberwithin{table}{section}

\title[On the topology of sums]
{On the topology of sums in powers \\ of an algebraic number}
\author{Nikita Sidorov}
\address{
School of Mathematics, The University of Manchester,
Oxford Road, Manchester M13 9PL, United Kingdom. E-mail:
sidorov@manchester.ac.uk}
\author{Boris Solomyak}
\address{Box 354350, Department of Mathematics, University of
Washington, Seattle, WA~98195, USA. E-mail:
solomyak@math.washington.edu}
\date{\today}
\subjclass[2000]{Primary 11J17; Secondary 11K16, 11R06} \keywords{Algebraic number, Perron number, Salem number, power sum.}

\begin{document}

\baselineskip=17pt

\begin{abstract}
Let $1<q<2$ and
\[
\La(q)=\left\{\sum_{k=0}^n a_kq^k\mid a_k\in\{-1,0,1\},\ n\ge1\right\}.
\]
It is well known that if $q$ is not a root of a polynomial with coefficients $0,\pm1$, then $\La(q)$ is dense in $\BR$. We give several sufficient conditions for the denseness of $\La(q)$ when $q$ is a root of such a polynomial. In particular, we prove that if $q$ is not a Perron number or it has a conjugate
$\alpha$ such that $q|\alpha|<1$, then $\La(q)$ is dense in $\BR$.
\end{abstract}

\maketitle

\section{Introduction and auxiliary results}

Let $q\in (1,2)$ and put
\[
\La_n(q)=\left\{\sum_{k=0}^{n} a_kq^k\mid a_k\in\{-1,0,1\}\right\},
\]
and $\La(q)=\bigcup_{n\ge1}\La_n(q)$. (It is obvious that the sets $\La_n(q)$ are nested.) The question we want to address is the topological structure of $\La(q)$. Is it dense? discrete? mixed?

The first important result has been obtained by A.~Garsia \cite{Ga}: if $q$ is a Pisot number (an algebraic integer greater than~1 whose conjugates are less than~1 in modulus), then $\La(q)$ is uniformly discrete. On the other hand, if $q$ does not satisfy an algebraic equation with coefficients $0,\pm1$, then it
is a simple consequence of the pigeonhole principle that $0$ is a limit
point of $\La(q)$ and thus, it is dense -- see below.

Surprisingly little is known about the case when $q$ is a root of a polynomial with coefficients $0,\pm1$. The most notable result is \cite[Theorem~I]{EK} in which the authors prove in particular that if $q<\frac{1+\sqrt5}2$ and $q$ is not Pisot, then $\La(q)$ has a finite accumulation point.

In this paper we study this case and give two sufficient conditions for $\La(q)$ to be dense. These conditions are rather general and cover a substantial subset of such $q$'s -- see Theorems~\ref{thm:1} and \ref{thm:main}.

Put
\[
Y_n(q)=\left\{\sum_{k=0}^{n} a_kq^k\mid a_k\in\{0,1\}\right\}
\]
and $Y(q)=\bigcup_{n\ge1}Y_n(q)$. The set $Y(q)$ is discrete and we can
write its elements in the ascending order:
\[
Y(q)=\{0=y_0(q) < y_1(q) < y_2(q)<\dots\}.
\]
Following \cite{EK}, we define
\[
l(q)=\liminf_{n\to\infty} (y_{n+1}(q)-y_n(q)).
\]

\begin{thm} (\cite{Dr})\label{thm:dr}
If\/ $0$ is a limit point of $\La(q)$, then $\La(q)$ is dense in $\BR$.
\end{thm}

It is obvious that $0$ is a limit point of $\La(q)$ if and only if $l(q)=0$. Hence follows

\begin{cor}The set $\La(q)$ is dense in $\BR$ if and only if $l(q)=0$.
\end{cor}

The purpose of this paper is to find some wide classes of algebraic $q$ for which $l(q)=0$.

Put for any $\be\in\BC$,
\[
Y_n(\be)=\left\{\sum_{k=0}^{n} a_k\be^k\mid a_k\in\{0,1\},\ 0\le k\le n
\right\}
\]
and $z_n(\be):=\# Y_n(\be)$. It is obvious that $z_n(\be)\le 2^{n+1}$.

In order to estimate $z_n(\beta)$ for $|\beta|>1$, it is
useful to consider the set
$$
A_\lam:= \left\{\sum_{k=0}^\infty a_k \lam^k\,|\ a_k \in \{0,1\},\ k\ge 0\right\},\ \ \mbox{where}\ \ \lam = \beta^{-1}.
$$
We have $|\lam|<1$, so the series converges for any choice of the coefficients $a_k \in \{0,1\}$. It is easy to see that the set
$A_\lam$ is
compact, being the image of the infinite product space $\{0,1\}^\infty$ under a continuous mapping.
It satisfies the set equation
$$
A_\lam = \lam A_\lam \cup (1+\lam A_\lam),
$$
and can be characterized as the unique compact set with this property \cite{Hutch}.
It is thus the
attractor of the iterated function system $\{z\mapsto \lam z,\ z\mapsto \lam z+1\}$ in the complex plane, see \cite{Hutch} for details.

The sets $A_\lam$, with $|\lam|<1$, have been extensively studied in the ``fractal'' literature; see e.g.\
\cite{Bandt,BH,IJK,SolXu} and the book \cite[8.2]{Barn}.
Note that some of these sources are concerned with the sets
$$
\wtil{A}_\lam:= \left\{\sum_{k=0}^\infty a_k \lam^k\,|\ a_k \in \{-1,1\},\ k\ge 0\right\},
$$
however, it is clear that $A_\lam = T(\wtil{A}_\lam)$, where $T(z) = \frac{1}{2}(z + (1-\lam)^{-1})$, so all the results immediately transfer.

\begin{lemma} \label{lem:fract}
\begin{enumerate}
\item If $\lam \in \BC$, with $|\lam|\in \bigl(\frac12,1\bigr)$, then $z_n(\lam) = \#Y_n(\lam) \ge |\lam|^{-n-1}$ for all $n$.
\item
If $\lam \in \BC$, with $2^{-1/2} \le |\lam|< 1$, and $|\Re \lam| \le |\lam|^2-\frac12$, then  $z_n(\lam)\ge |\lam|^{-2(n+1)}$ for all $n$.
\end{enumerate}
\end{lemma}
\begin{proof}
By the definition of the set $A_\lam$, we have for all $n\ge 0$:
\begin{equation} \label{eq:fract1}
A_\lam = \bigcup_{z\in Y_n(\lam)} (z + \lam^{n+1} A_\lam).
\end{equation}

\smallskip\noindent
{\rm (i)}\ Suppose that the set $A_\lam$ is connected, and let $u,v\in A_\lam$ be such that $|u-v|=\diam(A_\lam)$.
Then there exists a ``chain'' of distinct subsets $A_j:=z_j+\lam^n A_\lam\subset A_\lam$, $j=1,\ldots,m$, with $z_j \in Y_n(\lam)$, such that
$u\in A_1, v\in A_m$ and $A_j\cap A_{j+1}\ne \es$ for all $j\le m-1$. Therefore,
\begin{align*}
\diam(A_\lam) &\le \sum_{j=1}^m \diam(A_j) = m\cdot \diam(\lam^{n+1} A_\lam) \\
&\le \# Y_n(\lam) |\lam|^{n+1}\diam(A_\lam),
\end{align*}
and the claim follows.
If, on the other hand, $A_\lam$ is disconnected, then $\lam A_\lam \cap (\lam A_\lam + 1) = \es$. This is a general principle
for attractors
of iterated function systems with two contracting maps, see \cite{Hata,BH} or \cite[Chapter~8.2]{Barn}.
Therefore, in this case $\lam$ is not a zero of a power series with coefficients $\{-1,0,1\}$, much less a polynomial,
hence $z_n(\lam) = 2^{n+1} > |\lam|^{-n-1}$ for all $n$.

\medskip\noindent
{\rm (ii)}\
By \cite[Prop.\,2.6~(i)]{SolXu}, in view of the above remark concerning $\wtil{A}_\lam$, we know
that $A_\lam$ has nonempty interior for all $\lam$ in the open unit disc, such that $0 \le |\Re \lam| \le |\lam|^2-0.5$.
Then we have from (\ref{eq:fract1}) for the Lebesgue measure $\Lk^2$:
$$
\Lk^2(A_\lam) \le \# Y_n(\lam) \Lk^2(\lam^{n+1} A_\lam) = z_n(\lam)\cdot |\lam|^{2(n+1)} \Lk^2(A_\lam),
$$
as desired.
\end{proof}

Note that the proof of Lemma~\ref{lem:fract} did not use the fact that $\lam$ is non-real. Hence we obtain the following result as a direct corollary:

\begin{lemma}\label{lem:geq}
If $q\in(1,2)$, then $z_n(\pm q)\ge Cq^n$ for some $C>0$.
\end{lemma}

\begin{rmks}
\begin{enumerate}
\item Lemma~\ref{lem:geq} for $+q$ was proved in \cite{EK}, using the fact that $y_{n+1}(q)-y_n(q)\le 1$ for all $n$ and any $q\in(1,2)$.
\item With a bit more work one can show that in the setting of Lemma~\ref{lem:fract}~(i) we have $z_n(\lam) \ge C_n |\lam|^{-n}$ for some $C_n \uparrow
\infty$, assuming that $\lam$ is non-real. However, it is not needed in this paper.

 \item It follows from the results of \cite{DK,KL} that for any $\varphi \ne 0,\pi$, the set $A_\lam$ has nonempty interior for $\lam = re^{i\varphi}$, with $r$ sufficiently close to 1,
but it seems difficult to apply them in the absence of quantitative estimates.
\end{enumerate}
\end{rmks}

\begin{lemma}\label{lem2}
If $\be\in\BC\setminus\{0\}$, then $z_n(\be)=z_n(1/\be)$.
\end{lemma}
\begin{proof} Define $\phi:Y_n(\be)\to Y_n(1/\be)$ as follows:
\[
\phi\left(\sum_{k=0}^n a_k\be^k\right)=\sum_{k=0}^n a_{n-k}(1/\be)^k.
\]
A relation $\sum_{k=0}^n a_k\be^k = \sum_{k=0}^n b_k\be^k$ is equivalent to
$\sum_{k=0}^n a_k\be^{k-n} = \sum_{k=0}^n b_k\be^{k-n}$, which is in turn
equivalent to $\phi\left(\sum_{k=0}^n
a_k\be^k\right)=\phi\left(\sum_{k=0}^n b_k\be^k\right)$. Thus, $\phi$ is a
bijection.
\end{proof}

\begin{lemma}\label{lem:lq=0}
Let $q\in(1,2)$; if $z_n(q)\gg q^n$ (i.e., $\limsup_{n\to\infty}
q^{-n}z_n(q)=+\infty$), then $l(q)=0$.
\end{lemma}
\begin{proof} Since $\sum_{k=0}^n a_kq^k < q^{n+1}/(q-1)$, the result
follows immediately from the pigeonhole principle.
\end{proof}

Consequently, if $q$ is not a root of a polynomial with coefficients
$0,\pm1$, then $z_n(q)=2^{n+1}$, and $l(q)=0$ (which is well known, of course --
see, e.g., \cite{Dr}). If $q$ is such a root, it is obvious that $z_n(q)\ll
2^n$, and the problem becomes non-trivial. It is generally believed that
$l(q)=0$ unless $q$ is Pisot, but this is probably a very tough conjecture.

\section{Main results}

We need some preliminaries. Put
\[
L(q)=\limsup_{n\to\infty} (y_{n+1}(q)-y_n(q)).
\]
Note that $L(q)=0$ is equivalent to $y_{n+1}(q)-y_n(q)\to0$ as $n\to\infty$.
This condition was studied in the seminal paper \cite{EK}; in particular, it
was shown that if $q<2^{1/4}\approx1.18921$ and $q$ is not equal to the
square root of the second Pisot number $\approx 1.17485$, then $L(q)=0$\footnote{V.~Komornik has recently shown \cite{Kom-pers} that the second condition can be removed, so $L(q)=0$ if $q<2^{1/4}$.}. It was also shown in the same paper that $L(\sqrt2)=0$.

It is worth noting that the two conditions $l(q)=0$ and $L(q)=0$ are, generally speaking, very different in nature; for instance, as we know, $l(q)=0$ for all transcendental $q$, whereas $L(q)=1$ for all $q\ge\frac{1+\sqrt5}2$ (see, e.g.,  \cite{EJK}) and no $q\in\bigl(\sqrt2,\frac{1+\sqrt5}2\bigr)$ is known for which $L(q)=0$.

Throughout this section we assume that $q\in(1,2)$ is a root of a polynomial
with coefficients $0,\pm1$. It is easy to show that in this case any conjugate of $q$ is less than~2 in modulus.

Finally, recall that an algebraic integer $q>1$ is called a {\em
Perron number} if each of its conjugates is less than $q$ in modulus.

\begin{thm}\label{thm:1}
If $q\in (1,2)$ is not a Perron number, then $l(q)=0$. If, in addition, $q<\sqrt2$ and $-q$ is not a conjugate of $q$, then $L(q)=0$.
\end{thm}
\begin{proof} We first prove $l(q)=0$. We have three cases.

\smallskip\noindent
\textbf{Case 1.} $q$ has a real conjugate $p$ and $q<|p|$. Since $p$ is an
algebraic conjugate of $q$, it follows from the Galois theory that the map $\psi:Y_n(q)\to
Y_n(p)$ given by $\psi\left(\sum_{i=0}^n a_iq^i\right)=\sum_{i=0}^n a_i
p^i$, is a bijection. Hence $z_n(q)=z_n(p)\ge C|p|^n$ by Lemma~\ref{lem:geq}
and $z_n(q)\gg q^n$. Now the claim follows from Lemma~\ref{lem:lq=0}.

\smallskip\noindent
\textbf{Case 2.} $q$ has a complex non-real conjugate $p$ and $q<|p|$. This
case is similar to Case~1: $z_n(q)=z_n(p)\ge C|p|^n$ by
Lemma~\ref{lem:fract}~(i) and $z_n(q)\gg q^n$.

\smallskip\noindent
\textbf{Case 3.} $q$ has a conjugate $p$ and $q=|p|$. Let $f$ denote the minimal polynomial for $q$. Then we have $f(x)=g(x^m)$ for some $m\ge2$ by \cite{boyd}.
Put $\be=q^m$. We have
\begin{align*}
Y_{mk}(q)&=\{a_0+a_1\be^{\frac1m}+a_2\be^{\frac2m}+\dots +a_{mk}\be^n\mid
a_i\in\{0,1\}\}\\
&=\Bigl\{A_1+\be^{\frac1m}A_2+\be^{\frac2m}A_3 +\dots
+ \be^{\frac{m-1}m}A_m : \\
 & \ \ \ \ \ \ A_1\in Y_k(\be), A_i\in Y_{k-1}(\be), \ 2\le i\le m\Bigr\}.
\end{align*}
Observe that any relation of the form \[
A_1+\be^{\frac1m}A_2+\dots
+ \be^{\frac{m-1}m}A_m=A'_1+\be^{\frac1m}A'_2+\dots
+ \be^{\frac{m-1}m}A'_m
 \]
implies $A_1=A_1',\dots, A_m=A_m'$. Indeed,
if $q$ satisfies an equation $ B_1 + q B_2 + ... + q^{m-1} B_m = 0 $
with $B_i\in \BZ[q^m]$, then $qe^{2\pi i j/m}$ satisfies the same
equation for $j=1,\ldots,m-1$, hence
$B_i=0$ for all $i$. Thus, $z_{mk}\bigl(\be^{\frac1m}\bigr)=z_k(\be)\cdot (z_{k-1}(\be))^{m-1}$.

Now, if $q\ge2^{\frac1m}$, then $\be \ge 2$, so $z_k(\be)=2^{k+1}$, and we obtain from the above argument that for $n=mk$ we have
$z_n(q)\ge C2^n\gg q^n$. Otherwise $z_n(q)\ge z_n(\be) \ge C\be^n\gg q^n$. Hence
by Lemma~\ref{lem:lq=0}, $l(q)=0$.

Let us now prove the second part of the theorem. Suppose $q<\sqrt2$ is not Perron and $-q$ is not its conjugate; then
$q$ has a conjugate $\alpha \ne -q$, with $|\alpha| \ge q$. Thus,
$q^2$ has a conjugate $\alpha^2$, and $|\alpha|^2\ge q^2$ with $\alpha^2\neq q^2$. If $|\alpha|>\sqrt2$, then $\alpha^2$ (and, consequently, $q^2$) is not a root of $-1,0,1$ polynomial. Otherwise, we can apply the first part of this theorem to $q^2$. In either case, $l(q^2)=0$, whence by \cite[Theorem~5]{EJK}, $L(q)=0$.
\end{proof}

\begin{rmk} Stankov \cite{St} has proved a similar result for the following set:
\begin{equation}\label{eq:Aq}
\mathcal A(q)=\left\{\sum_{k=0}^n a_kq^k\mid a_k\in\{-1,1\},\ n\ge1\right\}.
\end{equation}
More precisely, he has shown that if $\mathcal A(q)$ is discrete, then all {\em real} conjugates of $q$ are of modulus strictly less than $q$.
\end{rmk}

\begin{cor}If $q\in(1,2)$ is the square root of a Pisot number and not itself Pisot, then $l(q)=0$.
\end{cor}
\begin{proof}If $q=\sqrt\be$ and $\be$ is Pisot, then either $-q$ is a conjugate of $q$ or $q$ is Pisot.
\end{proof}

\begin{thm}\label{thm:main}
\begin{enumerate}
\item
Suppose $q\in (1,2)$ has a conjugate $\alpha$ such that
$|\alpha| q < 1$. Then $l(q)=0$ and, consequently, $\La(q)$ is dense in
$\BR$.
\item Suppose $q\in (1,2)$ has a non-real conjugate $\alpha$ such that $|\alpha| q=1$. Then $l(q)=0$.
\end{enumerate}

\smallskip\noindent If, in addition, $q<\sqrt2$ in either case, then $L(q)=0$.
\end{thm}
\begin{proof}
(i)
As above, we have $z_n(q)=z_n(\alpha)$. On the other hand, by
Lemma~\ref{lem2},
$z_n(\alpha)=z_n(1/\alpha)$, and by Lemmas~\ref{lem:geq} and
\ref{lem:fract},
$z_n(1/\alpha)\ge C\cdot (|1/\alpha|)^n$. Hence $z_n(q)\ge C\cdot
(|1/\alpha|)^n\gg q^n$, in view of $|\alpha q|<1$. Hence by
Lemma~\ref{lem:lq=0}, $l(q)=0$.

If $q<\sqrt2$, then $q^2$ has a conjugate $\alpha^2$, and $q^2|\alpha|^2<1$.
Hence $l(q^2)=0$, whence $L(q)=0$.

\smallskip \noindent (ii) Denote $\alpha_1 = q, \alpha_2 = \alpha$, and
$\alpha_3 = \overline{\alpha}$. Since $|\alpha| q=1$ and $\alpha$ is non-real,
we have three conjugates satisfying $\alpha_1^2\alpha_2\alpha_3=1$. Smyth \cite[Lemma 1]{Smy} characterizes such situations, but it is
easier for us to proceed directly. The Galois group of the minimal polynomial for $q$ is
transitive, so there is an automorphism of the Galois group mapping $\alpha_1$ to $\alpha_2$. We obtain that $\alpha_2^2 \alpha_i \alpha_j =1$ for some
distinct conjugates $\alpha_i$ and $\alpha_j$ of $\alpha_1$.
But this implies
$\max\{|\alpha_i|,|\alpha_j|\} \ge \alpha_1=q$, hence $q$ is not a Perron number, and $l(q)=0$ by Theorem~\ref{thm:1}.

If $q<\sqrt2$, then $q^2|\alpha^2|=1$, and the first part of (ii) applies to $q^2$, unless $\alpha^2\in\BR$. If this is the case, then $\alpha=\pm i/q$, whence the minimal polynomial for $q$ contains only powers divisible by 4. Hence the minimal polynomial for $q^2$ contains only even powers, which implies that $-q^2$ is conjugate to $q^2$, whence $q^2$ is not Perron, and $l(q^2)=0$.
\end{proof}

\begin{rmk}
If $|\alpha| q=1$ and $\alpha$ is real, we do not know if $l(q)=0$. In fact, this includes the interesting (and probably, difficult) case of Salem numbers\footnote{Recall that an algebraic  number $q>1$ is called a {\em Salem number} if all its conjugates have absolute value no greater than~1, and at least one has absolute value exactly~1.}.
\end{rmk}

\begin{Def}We say that an algebraic integer $q>1$ is {\em anti-Pisot} if it has only one conjugate less than~1 in modulus and at least one conjugate greater than~1 in modulus other than $q$ itself.
\end{Def}

\begin{cor}\label{cor:anti}
If $q\in(1,2)$ is anti-Pisot and also a root of a polynomial with the coefficients in $\{-1,0,1\}$, then $l(q)=0$.
\end{cor}
\begin{proof}Let $\alpha=\alpha_1,\alpha_2,\dots,\alpha_{k-1},q$ be all the
conjugates of $q$. We have $\left|\prod_{j=1}^{k-1}\alpha_j\right|\cdot
q=1$, because $q$ satisfies an algebraic equation with coefficients
$0,\pm1$, whence its minimal polynomial must have a constant term $\pm1$.

Suppose $|\alpha|<1$; then it is clear than $\alpha\in\BR$ (since it is unique). If $|\alpha_2|>1$ and $|\alpha_j|\ge1$ for $j=3,\dots,k-1$, then it is obvious that $|\alpha|q\le|\alpha_2|^{-1}<1$, i.e., the condition of
Theorem~\ref{thm:main}~(i) is satisfied.
\end{proof}

\section{Examples}\label{sec:ex}

\begin{example} Let $q\approx1.22074$ be the positive root of $x^4=x+1$. Then $q$ has a single conjugate $\alpha\approx-0.72449$ inside the open unit disc and no conjugates of modulus~1, whence $q$ is anti-Pisot, and by Corollary~\ref{cor:anti}, $l(q)=0$. Furthermore, $q<\sqrt2$, whence $L(q)=0$ as well.

Note that $q>2^{1/4}$, so we cannot derive the latter claim immediately from \cite[Theorem~IV]{EK}.
\end{example}

\begin{example}An example of $q$ with a real conjugate $\alpha$ which is not anti-Pisot but still satisfies the condition of Theorem~\ref{thm:main}~(i), is the appropriate root of $x^5=x^4+x^2+x-1$. Here $q\approx1.52626$ and $\alpha\approx0.59509$.
\end{example}

\begin{example}For the equation $x^5=x^4-x^2+x+1$ we have $q\approx1.26278$
and $|\alpha|\approx0.74090$ so $|\alpha|q\approx0.93559$ (and
$\alpha\notin\BR$). By Theorem~\ref{thm:main}~(i), $L(q)=0$.
\end{example}

\begin{example} For the equation $x^8=x^7 + x^6 + x^5 - x^4-x^3-x^2+x-1$ we have $q\approx 1.52501$. Among its conjugates is
$\alpha\approx 0.3741 + 0.52404i$ with $|\alpha|\approx 0.64387 < 1/q = 0.65574$, so again $l(q)=0$ by Theorem~\ref{thm:main}~(i). Note that $q>\sqrt2$ so we cannot claim $L(q)=0$.
\end{example}

\begin{example}\label{ex:salem}
The following example illustrates Theorem~\ref{thm:main}~(ii). Let $q\approx1.19863$ be the largest root of $x^{12}=x^9+x^6+x^3-1$; then $\alpha=\zeta q^{-1}$ is a root of this equation as well, where $\zeta$ is any complex non-real cubic root of unity. Hence $q|\alpha|=1$, and Theorem~\ref{thm:main}~(ii) applies, i.e., $L(q)=0$. Note that $q=\sqrt[3]\beta$, where $\be$ is a quartic Salem number.
\end{example}

\begin{example} For the equation $x^{11} = x^{10} + x^9 - x^6 + x^4 - x^2 - 1$ we have $q\approx 1.5006$. Among its conjugates is
$\lam\approx 0.02625 + 0.7414i$. Theorem~\ref{thm:main} does not apply, but we can use Lemma~\ref{lem:fract}~(ii) to obtain
$$
z_n(q) = z_n(\lam) \ge |\lam|^{-2(n+1)} \approx  1.81696^{n+1} \gg q^n,
$$
which implies that $l(q)=0$. Note that Lemma~\ref{lem:fract}~(ii) applies, because $0.02625\approx \Re \lam < |\lam|^2-\frac12\approx 0.05037$.
\end{example}

\begin{example}Consider the equation
$x^{18}=-x^{16}+x^{14}+x^{11}+x^{10}+\dots+x+1$ (no powers missing between $x^{10}$ and 1). It has a root $q\approx 1.22289$, and the largest in modulus conjugates are $u,\overline u$ approximately equal to $-.03958\pm 1.3109i$. Then Theorem~\ref{thm:1} implies $L(q)=0$.

It is worth mentioning that there is another way to obtain this result. Consider $q^2$ and its conjugates $u^2, \overline u^2$. We claim that although $|u^2|<2$, $u^2$, and hence
$q^2$, is not a zero of a $-1,0,1$ polynomial (whence $l(q^2)=0$, which implies $L(q)=0$).

Indeed, if it were, then $q^{-2}, u^{-2}, (\overline u)^{-2}$ would also be zeros of such a polynomial. However, the product of these three numbers is
$\approx 0.226024$, so this is impossible, in view of the following

\medskip\noindent
{\bf Claim.} {\em Suppose $z_1, z_2, z_3$ are three different roots of a $-1,0,1$ polynomial. Then $|z_1 z_2 z_3| \ge 1/2\cdot(4/3)^{-3/2} =0.32476\dots$}

\smallskip
This claim is a slight generalization of \cite[Theorem~2]{borwein}, see \cite[Theorem 2.4]{shmerkin}.

\end{example}

\begin{example}\label{ex:bad}
Finally, an example of $q$ for which none of our criteria works is the real root of $x^5=x^4+x^3-x+1$. Here $q\approx1.54991$, and the other four
conjugates are non-real, with the moduli $\approx1.04492$ and $\approx0.76871$ respectively.

Another example is any Salem number $q\in(1,2)$, for instance $q\approx1.72208$ which is a root of $x^4=x^3+x^2+x-1$. (Which is of course none other than $\be$ from Example~\ref{ex:salem}.)
\end{example}

\section{Final remarks and open problems}

\subsection{} Our first remark concerns the case $q\in(m,m+1)$ with $m\ge2$.
Here the natural definition for $\La(q)$ is
\[
\La(q)=\left\{\sum_{k=0}^n a_kq^k\mid a_k\in\{-m,-m+1,\dots,m-1,m\},\
n\ge1\right\}.
\]
Theorem~\ref{thm:main} holds for this case, provided $\alpha\in\BR$ (and so does Case~1 of Theorem~\ref{thm:1})--- the proof is essentially the same.
The case of non-real $\alpha$ is less straightforward, since there is no ready-to-apply complex machinery for $m\ge2$.
(Basically, we need that if $\alpha$ is a zero of a polynomial with coefficients in $\{-m,\ldots,m\}$, then the attractor of the iterated function system
$\{\alpha z + j\}_{j=0}^m$ in the complex plane is connected. This can be verified for $m=2,3$ but we do not know if this is true in general.)
Note also that an analogue of Theorem~\ref{thm:dr} for $m\ge2$ has been proved in \cite{DrMc}.

\subsection{} We do not know whether the extra condition that $-q$ is not a conjugate of $q$ is really necessary in the second claim of Theorem~\ref{thm:1}. In particular, is it true that $L(\sqrt\varphi)=0$ if $\varphi$ is the golden ratio?

\subsection{} In \cite[Proposition~1.2]{PS} it is shown that if $q<\sqrt2$ and $q^2$ is not a root of a polynomial with coefficients $0,\pm1$, then the set $\mathcal A(q)$ given by (\ref{eq:Aq}) is dense in $\BR$. In fact, what the authors use in their proof is the condition $l(q^2)=0$. Consequently, Theorems~\ref{thm:1} and \ref{thm:main} provide sufficient conditions for $\mathcal A(q)$ to be dense in case when
$q^2$ does satisfy an algebraic equation with coefficients $0,\pm1$.

\subsection{} Is $l(q)=0$ for $q$ in Example~\ref{ex:bad} and suchlike?

\subsection{} All our criteria yield that $l(q)=0$ implies $L(q)=0$ for $q<\sqrt2$. Is this really the case?

\medskip\noindent\textbf{Acknowledgement.} The authors are indebted to the Max-Planck Institute where a significant part of the work was done in the summer of 2009. We are also grateful to Martijn de Vries for indicating the papers \cite{Dr, DrMc}.
The research of Solomyak was supported in part by NSF grant DMS-0654408.

\medskip\noindent\textbf{Added in proof.} In the recent paper by Sh.~Akiyama and V.~Komornik \cite{AK} several results mentioned in the introductory part of the present paper have been significantly improved, namely:

\begin{enumerate}
\item If $q\in(1,\sqrt2]$ is non-Pisot, then $\ell(q)=0$ and $\mathcal A(q)$ is dense in $\mathbb R$;
\item If $q\in(\sqrt2,2)$ is non-Pisot, then $\Lambda(q)$ has a finite accumulation point;
\item For $q\in(1,2^{1/3}]$ we have $L(q)=0$.
\end{enumerate}

\end{document}